\documentclass[11pt]{amsart}

\usepackage{amssymb}
\usepackage{latexsym}
\usepackage{amsmath}
\usepackage{amsthm}
\usepackage{amsfonts}
\usepackage{color}
\usepackage{pdfsync}
\usepackage[usenames,dvipsnames]{pstricks}
\usepackage{epsfig}
\usepackage{pst-grad} 
\usepackage{pst-plot} 

\newcommand{\R}{\mathbb{R}}

\newcommand{\I}{\mathcal{I}}
\newcommand{\U}{\mathcal{U}}

\theoremstyle{plain}
\newtheorem{defi}{Definition}[section]
\newtheorem{prop}[defi]{Proposition}
\newtheorem{teo}[defi]{Theorem}
\newtheorem{cor}[defi]{Corollary}
\newtheorem{lema}[defi]{Lemma}

\newtheorem{remark}[defi]{Remark}

\theoremstyle{definition}

\theoremstyle{remark}

\numberwithin{equation}{section}

\begin{document}

\title[]{Uniform Equicontinuity for  a family of Zero Order operators approaching the fractional Laplacian.}

\author[]{Patricio Felmer}
\address{
Patricio Felmer - 
Departamento de Ingenier\'\i a Matem\'atica and CMM (UMI 2807 CNRS), Universidad de Chile, Casilla 170 Correo 3, Santiago, CHILE. (pfelmer@dim.uchile.cl)
}

\author[]{Erwin Topp}
\address{
Erwin Topp - 
Departamento de Ingenier\'\i a Matem\'atica (UMI 2807 CNRS), Universidad de Chile, Casilla 170, Correo 3, Santiago, CHILE; 
and Laboratoire de Math\'ematiques et Physique Th\'eorique (CNRS UMR 6083), F\'ederation Denis Poisson, 
Universit\'e Fran\c{c}ois Rabelais, Parc de Grandmont, 37200, Tours, FRANCE  (etopp@dim.uchile.cl).}

\date{\today}

\begin{abstract} In this paper we consider a smooth bounded  domain $\Omega \subset \R^N$  and a parametric family of radially symmetric  kernels 
$K_\epsilon: \R^N \to \R_+$ such that, for each $\epsilon \in (0,1)$, its $L^1-$norm is finite but it blows up as $\epsilon \to 0$. 
Our aim is to establish an $\epsilon$ independent modulus of continuity in ${\Omega}$,  for the solution $u_\epsilon$ 
of the homogeneous Dirichlet problem 
\begin{equation*}
\left \{ \begin{array}{rcll} - \I_\epsilon [u] &=& f & \mbox{in} \ \Omega. \\
u &=& 0 & \mbox{in} \ \Omega^c, \end{array} \right .
\end{equation*}
where $f \in C(\bar{\Omega})$ and the operator $\I_\epsilon$ has the form
\begin{equation*}
\I_\epsilon[u](x) = \frac12\int \limits_{\R^N} [u(x + z) + u(x - z) - 2u(x)]K_\epsilon(z)dz
\end{equation*}
and it approaches the fractional Laplacian as $\epsilon\to 0$.  The modulus of continuity is obtained combining the comparison principle with the translation invariance of $\I_\epsilon$, constructing suitable  barriers  that allow to manage the discontinuities that the solution $u_\epsilon$ may have on $\partial \Omega$. Extensions of 
this result to fully non-linear elliptic and  parabolic operators  are also discussed.
\end{abstract}

\maketitle

\section{Introduction.}

Let $\Omega \subset \R^N$ be a bounded open domain with $C^2$ boundary, $f \in C(\bar{\Omega})$ and $\epsilon \in (0,1)$. 
In this paper we are concerned on study of the Dirichlet problem
\begin{eqnarray}
\label{eq}- \I_\epsilon [u] & = & f  \quad \mbox{in} \ \Omega, \\
\label{exteriordata} u & = & 0 \quad \mbox{in} \ \Omega^c,
\end{eqnarray}
where $\I_\epsilon$ is a nonlocal operator approaching the fractional Laplacian as $\epsilon$ approaches $0$. 
We  focus our attention on $\I_\epsilon$ with the form
\begin{eqnarray}\label{Iepsilonorden1}
\I_\epsilon [u](x) := \int_{\R^N} [u(x + z) - u(x)] K_\epsilon(z)dz, 
\end{eqnarray}
where, for $\sigma \in (0,1)$ fixed, $K_\epsilon$ is defined as
\begin{eqnarray*}
K_\epsilon(z) := \frac{1}{\epsilon^{N + 2\sigma} + |z|^{N + 2\sigma}} = \epsilon^{-(N + 2\sigma)} K_1(z/\epsilon). 
\end{eqnarray*}

Notice that for each $\epsilon \in (0,1)$, $K_\epsilon$ is integrable in $\R^N$ with $L^1$ norm
equal to $C \epsilon^{-2\sigma}$, where $C>0$ is a  constant depending only on $N$ and $\sigma$. 
We point out that operators with kernel in $L^1$, like  $\I_\epsilon$, are known  in the literature as \textsl{zero order nonlocal operators}.

Operator $\I_\epsilon$ is a particular case of a broad class of nonlocal elliptic operators. 
In fact, given a positive measure $\mu$ satisfying the \textsl{L\'evy condition}
\begin{equation*}
\int_{\R^N} \min \{1, |z|^2\} \mu(dz) < \infty,
\end{equation*}
and, for each $x \in \R^N$ and $u: \R^N \to \R$ bounded and sufficiently smooth at $x$, the operator $\I_\mu[u](x)$ defined as
\begin{equation}\label{operadorgeneral}
\I_\mu[u](x) =  \int_{\R^N} [u(x + z) - u(x) - \mathbf{1}_{B_1(0)}(z) \langle Du(x), z \rangle]\mu(dz), 
\end{equation}
has been a subject of study in a huge variety of contexts such as potential theory (\cite{Landkof}), probability (\cite{Blumenthal, Sato})
and analysis (\cite{Sayah1, Sayah2, Barles-Imbert, Caffarelli-Silvestre1, Caffarelli-Silvestre2}). 
An interesting point of view of our problem 
comes from probability, since~\eqref{operadorgeneral} represents the infinitesimal generator 
of a \textsl{jump L\'evy process}, see Sato~\cite{Sato}. 
In our setting, the finiteness of the measure is associated with the so-called \textsl{Compound Poisson Process}. 
Dirichlet problems with the form of~\eqref{eq}-\eqref{exteriordata}
arise in the context of exit time problems with trajectories driven by the jump L\'evy process defined by $K_\epsilon(z)dz$,
and the solution $u_\epsilon$ represents the expected value of the associated cost functional, see~\cite{Oksendal}.

We may start our discussion with a natural notion of solution to our problem ~\eqref{eq}-\eqref{exteriordata}: we say that a {bounded} function $u: \R^N \to \R$, continuous in $\Omega$, is a   {solution} 
of~\eqref{eq}-\eqref{exteriordata} if it satisfies \eqref{eq}  pointwise in $\Omega$ and $u = 0$ on $\Omega^c$.
As we see in Section \S 2, this problem has a unique solution, more interestingly, through an example we will see that such a solution 
may not be  continuous in $\R^N$, since a discontinuity may appear on the boundary of $\Omega$. See Remark~\ref{rmku>0borde}.

This situation is in great contrast with  the limit case $\epsilon = 0$, where the kernel becomes $K(z)dz = |z|^{-(N + 2\sigma)}dz$ 
and the associated nonlocal operator is the \textsl{fractional Laplacian of order $2\sigma$}, denoted by $-(-\Delta)^\sigma$, see~\cite{Hitch}.
In this case, the corresponding Dirichlet problem becomes
\begin{eqnarray}\label{eqfractional}
\left \{ \begin{array}{rll} C_{N, \sigma} (-\Delta)^\sigma v &= f  \quad & \mbox{in} \ \Omega, \\ v &= 0 \quad & \mbox{in} \ \Omega^c,\end{array} \right .
\end{eqnarray}
where $C_{N,\sigma} > 0$ is a normalizing constant. 
In the context of the viscosity theory for nonlocal 
equations (see~\cite{Barles-Imbert,Sayah1,Sayah2}), Barles, Chasseigne and Imbert~\cite{Barles-Chasseigne-Imbert} addressed
a large variety of nonlocal elliptic problems including~\eqref{eqfractional}. In that paper, the authors proved
the existence and uniqueness of a viscosity solution $v \in C(\bar{\Omega})$ 
of~\eqref{eqfractional} satisfying $v = 0$ on $\partial \Omega$ that is, consequently,  continuous when we regard it as a
function on $\R^N$.  This result is accomplished by the use of a nonlocal version of  the notion of  viscosity solution with \textsl{generalized boundary conditions}, 
see~\cite{usersguide, Barles-Burdeau, Barles-Rouy} for an introduction of this notion in the context of second-order equations.

%
%

Additionally, fractional problems like~\eqref{eqfractional} enjoy a \textsl{regularizing effect} as in the classical second-order case. Roughly 
speaking, for a right-hand side which is merely bounded, the solution $v$ of~\eqref{eqfractional} is locally H\"older continuous in $\Omega$, 
see~\cite{Silvestre}. In fact, we should mention here that interior H\"older regularity for more general fractional problems 
(for which~\eqref{eqfractional} is a particular case) has been addressed by many authors, 
see for instance~\cite{Barles-Chasseigne-Imbert, Barles-Chasseigne-Imbert-Ciomaga, Bass-Kassman-Holder, 
Caffarelli-Silvestre1, Caffarelli-Silvestre2, Cardaliaguet-Cannarsa, Silvestre} and the classical book of Landkof~\cite{Landkof}, for 
 a non-exhaustive list of references. The interior  H\"older regularity is accomplished by well established elliptic techniques as the 
Harnack's inequality (\cite{Caffarelli-Silvestre1, Bass-Levin}) and the Ishii-Lions 
method (\cite{Barles-Chasseigne-Imbert, Ishii-Lions}). In both cases, 
 the nonintegrability of the kernel plays a key role. 
H\"older regularity for problems like~\eqref{eqfractional} can be extended up to the boundary, as it is proved by  Ros-Oton and Serra in \cite{Ros-Oton-Serra}, where a boundary Harnack's inequality is the key ingredient (see also~\cite{Bogdan}).
Naturally, as a byproduct of these regularity results, compactness properties are available for certain families of solutions 
of fractional equations. For instance, the family $\{ v_\eta \}$ of functions solving
\begin{equation*}
\left \{ \begin{array}{rll} C_{N, \sigma} (-\Delta)^\sigma v_\eta &= f_\eta  \quad & \mbox{in} \ \Omega \\ 
u_\eta &= 0, \quad & \mbox{in} \ \Omega^c, \end{array} \right . 
\end{equation*}
satisfies compactness properties when $\{ f_\eta \}$ is uniformly bounded in $L^\infty(\bar{\Omega})$.

For zero order problems, regularizing effects 
as arising  in fractional problems are no longer available (see~\cite{Chasseigne-Chaves-Rossi}). In fact, 
the finiteness of the kernel of zero order operators turns into degenerate ellipticity for which Ishii-Lions method cannot be applied.
Thus, ``regularity results'' for zero order problems like~\eqref{eq}-\eqref{exteriordata} are circumscribed to the heritage of 
the modulus of continuity of the right-hand side $f$ to the solution $u_\epsilon$ as it can be seen in~\cite{Chasseigne}. However, 
the modulus of continuity found  in~\cite{Chasseigne} depends strongly on the size  of the $L^1$ norm of $K_\epsilon$, which explodes as $\epsilon \to 0$.
A similar lack of stability as $\epsilon \to 0$ can be observed in the Harnack-type inequality results for nonlocal problems found 
by Coville in~\cite{Coville}. Hence, none of the mentioned tools are adequate for getting compactness for the family of solutions 
$\{u_\epsilon \}$ of problem~\eqref{eq}-\eqref{exteriordata},
which is a paradoxical situation since, in the limit case, the solutions actually get higher regularity and stronger compactness control 
on its behavior.

In view of the discussion given above, a natural mathematical question is if there exists a  uniform modulus of continuity  in $\Omega$, for the  family of  
solutions $\{u_\epsilon\}$ to \eqref{eq}-\eqref{exteriordata}, and consequently  compactness properties for it.  
In this direction, the main result of this paper is the following
\begin{teo}\label{teo1}
Let $\Omega \subset \R^N$ a bounded domain with $C^2$ boundary and $f \in C(\bar{\Omega})$. For $\epsilon \in (0,1)$, 
let $u_\epsilon$ be a solution to problem~\eqref{eq}-\eqref{exteriordata}. 
Then, there is a modulus of continuity $m$ depending only on $f$, such that
\begin{eqnarray*}
|u_\epsilon(x) - u_\epsilon(y)| \leq m(|x - y|), \quad \mbox{for } x,y \in \Omega.
\end{eqnarray*}
\end{teo}

The proof of this theorem  is obtained combining the translation invariance of $\I_\epsilon$ and comparison principle,
constructing suitable barriers to manage the discontinuities that $u_\epsilon$ may have on $\partial \Omega$ and to understand
how they evolve as $\epsilon$ approaches zero, see Proposition~\ref{propaprioriborde}.

As a consequence  of Theorem~\ref{teo1} we have the following corollary, that actually was our original motivation to study the problem. 
\begin{cor}\label{corolarioteo}
Let  $u_\epsilon$ be the solution to equation (\ref{eq}), with $f$ and $\Omega$ as in Theorem~\ref{teoborde}, and let $u$ be the solution of the equation (\ref{eqfractional}),
then $u_\epsilon \to u$ in $L^\infty(\bar{\Omega})$ as $\epsilon \to 0$.
\end{cor}

We mention here that the application of the 
\textsl{half-relaxed limits method} introduced by Barles and Perthame in~\cite{Barles-Perthame1} 
(see also \cite{Barles-Perthame, Bensaoud-Sayah, Barles-Imbert})
allows to obtain in a very direct way locally uniform convergence in $\Omega$ in the above corollary. 
At this point we emphasize on the main contribution of this paper, which is the analysis of the boundary behavior of the family
$\{ u_\epsilon \}$ of solutions to~\eqref{eq}-\eqref{exteriordata} 
coming from Theorem~\ref{teo1} and the subsequent \textsl{global} uniform convergence to the solution of~\eqref{eqfractional}. 

There are many possible extensions of Theorem \ref{teo1}, for example,
{it can be readily extended to problems with the form 
\begin{eqnarray*}
\left \{ \begin{array}{rll} - \I_\epsilon [u] &= f_\epsilon  \quad & \mbox{in} \ \Omega \\
u &= 0 \quad & \mbox{in} \ \Omega^c, \end{array} \right .
\end{eqnarray*}
with $\{f_\epsilon\} \subset  C(\bar{\Omega})$ having a common modulus of continuity independent of $\epsilon \in (0,1)$.} It can also be extended to
 fully nonlinear operators and to parabolic equations, as we discuss in   Section \S \ref{remarks}.  We could also consider different  families of approximating zero order operators, 
but we do not pursue this direction. 
There are many other interesting lines of research  that arises from this work. From the discussion given before 
Theorem~\ref{teo1}, questions arises with respect to Harnack type inequalities and its relation with regularity and compactness properties of solutions, when $\epsilon\to 0$. Regarding operators ${\mathcal I}_\mu$, where $\mu$ might be singular with respect to the Lebesgue measure, an interesting question that arises is if the main results of this article can be extended to this case.

The paper is organized as follows: In Section \S \ref{ppiosdelmaximo} we establish the notion of pointwise solution and the comparison principle. 
Important estimates for the discontinuity of the solution at the boundary are given in Section \S\ref{estimacionbordesection}, and
the boundary equicontinuity result is presented in Section \S \ref{seccionND}. The interior modulus of continuity is easily derived 
from the boundary equicontinuity, and therefore the proof of Theorem~\ref{teo1} is given in Section \S\ref{pruebamodulointerior}. 
Further related results are discussed in Section \S \ref{remarks}.

\subsection{Notation} For $x \in \R^N$ and $r>0$, we denote $B_r(x)$ the ball centered at $x$ with radius $r$ and simply $B_r$ if $x = 0$.
For a  set $U \subset \R^N$, we denote by $d_U(x)$ the signed distance to the boundary, this is $d_U(x) = \mathrm{dist}(x, \partial U)$, 
with $d_U(x) \leq 0$ if $x \in U^c$. Since many arguments in this paper concerns the set $\Omega$, we write $d_\Omega = d$. We also define
$$
\Omega_r = \{ x \in \Omega \ : \ d(x) < r \}
$$

Concerning the regularity of the boundary of $\Omega$, we assume it is at least $C^2$, so the distance function $d$ is a $C^2$ function in 
a neighborhood of $\partial \Omega$.
More precisely, there exists $\delta_0 > 0$ such that $x \mapsto d(x)$ is of class $C^2$ for $-\delta_0 < d(x) < \delta_0$. 

In our estimates we will denote by $c_i$ with $i=1, 2,...$ positive constants appearing in our proofs, depending only on $N, \sigma$ and $\Omega$. 
When necessary we will make explicit the dependence on the parameters. The index will be reinitiated in each proof.

\section{Notion of Solution and Comparison Principle.}
\label{ppiosdelmaximo}

In the introduction we defined a  notion  of solution to problem \eqref{eq}-\eqref{exteriordata}, which is very natural 
for zero order operators and allows us to understand the main features of the mathematical problem that we have at hand. However, this notion 
is not suitable for a neat statement of the comparison principle and it is not adequate to understand the limit as $\epsilon\to 0$. For this reason, from now on, we adopt another notion of solution which is more adequate, that is the notion of viscosity solution with generalized boundary condition defined by  Barles, Chasseigne and Imbert in~\cite{Barles-Chasseigne-Imbert}.

%
 
We remark that results  provided in this section are adequate for problems  slightly more general than 
our problem~\eqref{eq}-\eqref{exteriordata}. 
We will consider $J \in L^1(\R^N)$ a 
nonnegative function, and we define the nonlocal operator associated to $J$ as
\begin{equation}\label{IJ}
\I_J[u](x) = \int_{\R^N} [u(x + z) - u(x)] J(z)dz,
\end{equation}
for $u \in L^\infty(\R^N)$ and $x \in \R^N$, and a Dirichlet problem of the form
\begin{eqnarray}\label{eqbasica}
-\I_J[u] = f \quad & \mbox{in} \ \Omega,\\ \label{eqbasica1}
u=0\quad & \mbox{on} \ \Omega^c, 
\end{eqnarray}
with $f \in C(\bar{\Omega})$. 
Since we are interested in a Dirichlet problem for which the exterior data plays a role, 
we assume $J$ and $\Omega$ satisfy the condition
\begin{equation}\label{intJexterior}
\inf_{x \in \bar{\Omega}} \int_{\Omega^c - x} J(z)dz \geq \nu_0 > 0.
\end{equation}
Notice that problem~\eqref{eq}-\eqref{exteriordata} is a particular case of~\eqref{eqbasica}-\eqref{eqbasica1}.

In this situation, a bounded function $u: \R^N \to \R$, continuous in $\bar{\Omega}$ is a viscosity solution  with generalized boundary condition to problem \eqref{eqbasica}-\eqref{eqbasica1} if and only if
it satisfies
\begin{eqnarray}
\label{eqPtoFijo} -\I_J[u] & = f \quad & \mbox{on} \ \bar{\Omega}, \\
\label{exteriordataPtoFijo} u & = 0 \quad & \mbox{in} \ \bar{\Omega}^c.
\end{eqnarray}

The sufficient condition is direct from the definition and the necessary condition follows from the lemma:
\begin{lema}\label{lemahastaelborde}
Let $f \in C(\bar{\Omega})$ and let $u: \R^N \to \R$ be a function satisfying
\begin{equation}\label{ineqlemahastaelborde}
-\I_J[u](x) \leq f(x) \quad \mbox{for all} \ x \in \Omega, 
\end{equation}
where the above inequality is understood pointwise. Let $x_0 \in \partial \Omega$ and assume there exists a sequence $\{x_k\} \subset \Omega$ such that
\begin{eqnarray}\label{xktox0}
 x_k \to x_0, \quad u(x_k) \to u(x_0) 
\end{eqnarray}
and 
\begin{equation}\label{limsupprop}
\limsup_{k \to +\infty} u(x_k + z) \leq u(x_0 + z), \quad a.e.
\end{equation}

Then, $u$ satisfies~\eqref{ineqlemahastaelborde} at $x_0$.
\end{lema}

Here and in what follows  the considered measure  is the Lebesgue measure.

\medskip

\noindent
{\bf Proof.} Consider $\{x_k\} \subset \Omega$ as in~\eqref{xktox0}. Then, we can write
$$
\int_{\R^N} u(x_k + z)J(z)dz - u(x_k) \int_{\R^N}J(z)dz \geq -f(x_k).
$$
Hence, taking limsup in both sides of the last inequality, by~\eqref{xktox0} and the continuity of $f$, we arrive to
\begin{equation*}
\int_{\R^N} \limsup_{k \to \infty} u(x_k + z)J(z)dz - u(x_0) \int_{\R^N}J(z)dz \geq -f(x_0),
\end{equation*}
where the exchange of the integral and the limit is justified by Fatou's Lemma. Then, using~\eqref{limsupprop}, we conclude the result.
\qed

\medskip

We continue with our analysis  with an existence result for \eqref{eqbasica}-\eqref{eqbasica1}.
\begin{prop}\label{existencia}
Let $f \in C(\bar{\Omega})$. Then, there exists a unique bounded function $u: \R^N \to \R$, continuous in $\bar{\Omega}$, which is a viscosity solution with generalized boundary condition to problem \eqref{eqbasica}-\eqref{eqbasica1}.
\end{prop}
\noindent
{\bf Proof.} According with our discussion above, we need to find a solution to \eqref{eqPtoFijo}-\eqref{exteriordataPtoFijo}. Consider the map
$
T_a : C(\bar{\Omega}) \to C(\bar{\Omega})
$
defined as 
$$
T_a(u)(x) = u(x) - a\Big{(}||J||_{L^1(\R^N)}u(x) - \int_{\Omega - x} u(x + z)J(z)dz - f(x)\Big{)}.
$$
We observe that   $u\in C(\bar{\Omega})$ is a fixed points of $T_a$ if and only if  $u$ is a solution to problem~\eqref{eqPtoFijo}-\eqref{exteriordataPtoFijo}. Therefore, the aim 
is to prove that for certain $a > 0$ small enough, the map $T_a$ is a contraction in $C(\bar{\Omega})$. 
By~\eqref{intJexterior}, there exists $\varrho_0>$ such that \begin{equation*}
||J||_{L^1(\R^N)} - ||J||_{L^1(\Omega - x)} \geq \varrho_0,\quad \mbox{for each } x \in \bar{\Omega}. 
\end{equation*}
Let $0<  a < \min \{ \varrho_0^{-1}, ||J||_{L^1(\R^N)}^{-1} \}$ and consider  $u_1, u_2 \in C(\bar{\Omega})$. Then, for all $x \in \bar{\Omega}$ we have 
\begin{equation*}
\begin{split}
T_a(u_1)(x) - T_a(u_2)(x) \leq & \ \Big{(} 1 - a||J||_{L^1(\R^N)} + a \int \limits_{\Omega - x} J(z)dz \Big{)}||u_1 - u_2||_\infty \\ 
\leq & \ \Big{(} 1 - a \varrho_0 \Big{)}||u_1 - u_2||_\infty, 
\end{split}
\end{equation*}
concluding that
$$
||T_a(u_1) - T_a(u_2)||_\infty \leq (1 - a \varrho_0)||u_1 - u_2||_\infty,
$$
that is, $T_a$ is a contraction in $C(\bar{\Omega})$. From here existence and uniqueness follow.
\qed

\begin{remark}
We observe that $u:\R^N\to\R$, a viscosity solution with generalized boundary condition to problem \eqref{eqbasica}-\eqref{eqbasica1}, may be redefined on the boundary $\partial \Omega$ as $u=0$, to obtain a solution to 
\eqref{eqbasica}-\eqref{eqbasica1}
in the sense defined in the introduction.
\end{remark}
\begin{remark}\label{rmku>0borde}
Let $u$ be a solution of \eqref{eqbasica}-\eqref{eqbasica1}
in the sense defined in the introduction, with  $f \geq \varrho_0 > 0$. Our purpose is to show that $u$ has a discontinuity on the boundary of $\Omega$. 
Let us assume, for contradiction, that $u:\R^N\to\R$ is a continuous function. 
 
Then $u \geq 0$ in $\Omega$, otherwise
there exists $x_0 \in \Omega$ such that $u(x_0) = \min_{\bar{\Omega}} \{u\} < 0$ and evaluating the equation at $x_0$ we arrive to
\begin{equation*}
u(x_0) \int_{\Omega^c - x_0} J(z)dz \geq -\I_J[u](x_0) = f(x_0),
\end{equation*}
which is a contradiction to~\eqref{intJexterior}.
Then, from the equation, we have  for each $x \in  \Omega$ the inequality
\begin{equation*}
-\I_J[u](x) = f(x) > \varrho_0.
\end{equation*}
Since $u$ and $f$ are continuous and $u=0$ on $\partial \Omega$ and using  that  $u \geq 0$ in $\Omega$, we obtain  that, for each $x \in \partial \Omega$ 
\begin{equation*}
0 \geq -\int_{\Omega - x} (u(x + z)-u(x))J(z)dz = -\I_J[u](x) = f(x) > \varrho_0,
\end{equation*}
which is a contradiction. Thus, $u > 0$ on $\partial \Omega$ which implies that $u$ is discontinuous on $\partial\Omega$.
\end{remark}

In what follows we prove that a solution  to~\eqref{eqbasica}-\eqref{eqbasica1},  in the sense defined in the introduction, can be extended continuously to $\bar\Omega$. 
\begin{prop}\label{propexistenciapuntual}
Let $f \in C(\bar{\Omega})$. Let $v: \R^N \to \R$ in $L^\infty(\R^N) \cap C(\Omega)$ be a  solution to~\eqref{eqbasica}-\eqref{eqbasica1}, 
and $u: \R^N \to \R$ in $C(\bar{\Omega})$ be the viscosity solution to~\eqref{eqbasica}-\eqref{eqbasica1} given by Proposition~\ref{existencia}.
Then, $u = v$ in $\Omega$.
\end{prop}

\noindent
{\bf \textit{Proof:}} By contradiction, assume the existence of a point in $\Omega$ where $u$ is different from $v$. Defining $w = u- v$, 
we will assume that
\begin{equation}\label{Mlemauv}
M := \sup \limits_{\Omega} \{ w \} > 0,
\end{equation}
since the case $\inf_{\Omega} \{ w \} < 0$ follows the same lines. Moreover, we assume that the supremum defining $M$ is not attained, since this 
is the most difficult scenario. 
Let $\eta > 0$ and let $x_\eta \in \Omega \setminus \Omega_\eta$ such that
\begin{equation*}
w(x_\eta) = \max_{\Omega \setminus \Omega_\eta} \{ w \},
\end{equation*}
where $\Omega_\eta$ was defined at the end of the introduction.
We clearly have $w(x_\eta) \to M$ as $\eta \to 0$ and since we assume $M$ is not attained, then $x_\eta \to \partial \Omega$ as $\eta \to 0$.
Now, using the equations for $u$ and $v$ at $x_\eta \in \Omega$, we can write
\begin{equation*}
-\int_{\Omega - x_\eta} [w(x_\eta + z) - w(x_\eta)] J(z)dz + w(x_\eta) \int_{\Omega^c - x_\eta} J(z)dz \leq 0,
\end{equation*}
and by~\eqref{intJexterior} and the fact that $w(x_\eta) \to M$ as $\eta \to 0$, we have
\begin{equation}\label{ineqlemauv}
-\int_{\Omega - x_\eta} [w(x_\eta + z) - w(x_\eta)] J(z)dz + \nu_0 M - o_\eta(1) \leq 0,
\end{equation}
where $o_\eta(1) \to 0$ as $\eta \to 0$. But writing
\begin{equation*}
\begin{split}
\int_{\Omega - x_\eta} [w(x_\eta + z) - w(x_\eta)] J(z)dz = & \int_{\Omega \setminus \Omega_\eta - x_\eta} [w(x_\eta + z) - w(x_\eta)] J(z)dz \\
& \ + \int_{\Omega_\eta - x_\eta} [w(x_\eta + z) - w(x_\eta)] J(z)dz,
\end{split}
\end{equation*}
by the boundedness of $w$ and the integrability of $J$, the second integral term in the right-hand side of the last equality is $o_\eta(1)$,
meanwhile, using the definition of $x_\eta$ we have the first integral is nonpositive. Thus, we conclude
\begin{equation*}
\int_{\Omega - x_\eta} [w(x_\eta + z) - w(x_\eta)] J(z)dz \leq o_\eta(1), 
\end{equation*}
and replacing this into~\eqref{ineqlemauv}, we arrive to
\begin{equation*}
\nu_0 M - o_\eta(1) \leq 0.
\end{equation*}
By making   $\eta \to 0$, we see that this contradicts \eqref{Mlemauv}, since $\nu_0 > 0.$ 
\qed

\medskip

As a consequence of the last proposition, we have the following
\begin{cor}\label{corexistence}
Let $f \in C(\bar{\Omega})$. Then, there exists a unique solution $v \in L^\infty(\R^N) \cap C(\Omega)$ 
to problem~\eqref{eqbasica}-\eqref{eqbasica1} in the sense defined in the introduction.
Moreover, $v$ is uniformly continuous in $\Omega$ and its unique continuous extension to $\bar{\Omega}$ 
coincides with the unique viscosity solution to~\eqref{eqbasica}-\eqref{eqbasica1}. 
\end{cor}

The main tool in this paper is the comparison principle, and here the so-called \textsl{strong} 
comparison principle is the appropriate version to deal with  
discontinuities at the boundary. 
\begin{prop}\label{comparison}{\bf (Comparison Principle)}
Assume $f \in L^\infty(\bar{\Omega})$. Let $u, v \in \R^N \to \R$ be bounded, upper and lower semicontinuous functions on $\bar{\Omega}$, respectively. 
Assume $u$ and $v$ satisfy
\begin{equation}\label{ineqcomparison}
-\I_J[u] \leq f \quad \mbox{and} \quad -\I_J[v] \geq f, \quad \mbox{on} \ \bar{\Omega}.
\end{equation}
If $u \leq v$ in $\bar{\Omega}^c$, then $u \leq v$ in $\bar{\Omega}$.
\end{prop}

\noindent
{\bf \textit{Proof:}} Assume by contradiction that there exists $x_0 \in \bar{\Omega}$ such that 
$$
(u - v)(x_0) = \max_{x \in \bar{\Omega}} \{ u- v \} > 0.
$$ 
Evaluating  inequalities in~\eqref{ineqcomparison} at $x_0$ and substracting them, denoting $w = u - v$, we arrive to
\begin{equation*}
- \int_{\Omega - x_0} [w(x_0 + z) - w(x_0)] J(z)dz - \int_{\Omega^c - x_0} [w(x_0 + z) - w(x_0)] J(z)dz \leq 0,
\end{equation*}
and therefore, using that $x_0$ is a maximum point for $w$ in $\Omega$ and that $w \leq 0$ in $\Omega^c$, we can write
\begin{equation*}
w(x_0) \int_{\Omega^c - x_0} J(z)dz \leq 0, 
\end{equation*}
and using~\eqref{intJexterior} we arrive to a contradiction with the fact that $w(x_0) > 0$. 
\qed

\medskip

As a  first consequence of this comparison principle, we obtain an a priori $L^\infty(\bar{\Omega})$ estimate 
for the solutions $u_\epsilon$ of~\eqref{eq}-\eqref{exteriordata}, independent of $\epsilon$.

\begin{prop}\label{cotauniformesolucion}
Let $\epsilon \in (0,1)$, $f \in C(\bar{\Omega})$ and $u_\epsilon$ be the viscosity solution of~\eqref{eq}-\eqref{exteriordata}. 
Then, there exists a constant ${C} > 0$ 
such that 
$$
||u_\epsilon||_{L^\infty(\bar{\Omega})} \leq C ||f||_\infty
$$
and this constant depends only on $\Omega$, $N$ and $\sigma$, but not on $\epsilon$, for $\epsilon\in (0,1)$.
\end{prop}
\noindent
{\bf Proof.} Consider the bounded function
$
\chi(x) = \mathbf{1}_{\bar{\Omega}}(x).
$
We clearly have that $\chi \in C(\bar{\Omega})$ and $\chi = 0$ in $\bar{\Omega}^c$. 
Denote $R = diam(\Omega) > 0$ and use the definition of the operator $\I_\epsilon$ to see
that for each $x \in \bar{\Omega}$ we have
\begin{equation*}
\begin{split}
-\I_\epsilon[\chi](x) = \int_{\Omega^c - x} K_\epsilon(z)dz \geq \int_{B_{R + 1}^c} \frac{dz}{2|z|^{N + 2\sigma}}
= \frac{\mathrm{Vol}(B_1)(R + 1)^{-2\sigma}}{2\sigma}.
\end{split}
\end{equation*}
Hence, denoting $C = (2\sigma)^{-1}\mathrm{Vol}(B_1)(R + 1)^{-2\sigma}$ and $\tilde{\chi} =C^{-1} ||f||_\infty \chi$, we may use the  comparison principle to conclude $u_\epsilon \leq C ||f||_\infty$
in $\bar{\Omega}$. A lower bound can be found in a similar way, concluding the result.
\qed


\section{Estimates of the Boundary Discontinuity.}
\label{estimacionbordesection}

The aim of this section is to estimate the discontinuity jump on $\partial \Omega$ of the solution $u_\epsilon$ of~\eqref{eq}-\eqref{exteriordata}. For this purpose,
 a flattening procedure on the boundary is required.

Recall that $\delta_0 > 0$ is such that the distance function to $\partial \Omega$ is smooth in $\Omega_{\delta_0}$, and for $x \in \Omega_{\delta_0}$
we denote $\hat{x}$ the unique point on $\partial \Omega$ such that $d_\Omega(x) = |x - \hat{x}|$.
We can fix $\delta_0$ small in order to have the existence of three constants $R_0, r_0, r_0' > 0$ depending only on the regularity 
of the boundary, satisfying the following properties:

\medskip
\noindent
{\bf (i)} For each $x \in \Omega_{\delta_0}$, there exists $\mathcal{N}_x \subset \partial (\Omega - x)$, a $\partial (\Omega - x)$-neighborhood 
of $\hat{x} - x$, which is the graph of a $C^2$ function $\varphi_x: B_{R_0} \subset \R^{N-1} \to \R^N$, that is, 
\begin{equation*}
(\xi', \varphi_x(\xi')) \in \mathcal{N}_x, \quad \mbox{for all} \ \xi' \in B_{R_0}. 
\end{equation*}

\medskip
\noindent
{\bf (ii)} If we define the function $\Phi_x$ as
\begin{equation}\label{Phi}
\Phi_x(\xi', s)  = (\xi', \varphi(\xi')) + (d(x) + s) \nu_{\xi'}, \quad (\xi',s) \in B_{R_0} \times (-R_0, R_0),
\end{equation}
where $\nu_{\xi'}$ is the unit inward normal to $\partial (\Omega - x)$ at $(\xi', \varphi_x(\xi'))$ and
denoting
$\mathcal{R}_x = \Phi_x(B_{R_0} \times (-R_0, R_0))$, then $\Phi_x : B_{R_0} \times (-R_0, R_0) \to \mathcal{R}_x$ 
is a $C^1$-diffeomorphism. Notice that $\Phi_x(0',0)$ is the origin and therefore $\mathcal{R}_x$ is an $\R^N$-neighborhood of the origin. 

\medskip
\noindent
{\bf (iii)} The constant $r_0 > 0$ is such that $B_{r_0} \subset \mathcal{R}_x$ for all $x \in \Omega_{\delta_0}$. 

\medskip
\noindent
{\bf (iv)} The constant $r_0' > 0$ is such that $\Phi_x (B_{r_0'} \times (-r_0', r_0')) \subset B_{r_0}$.

\medskip

We may assume $0 < r_0' \leq r_0 \leq \delta_0$.
In addition, by the smoothness of the boundary there exists
a constant $C_\Omega > 1$ such that
\begin{equation}\label{tildeK}
C_\Omega^{-1} K_\epsilon(\xi) \leq \tilde{K}_\epsilon(\xi) \leq C_\Omega K_\epsilon(\xi), \quad \xi \in \R^N,
\end{equation}
where $\tilde{K}_\epsilon(\xi) = |\mathrm{Det}(D\Phi_x(\xi))|K_\epsilon(\Phi_x(\xi))|$.

The following Lemma is the key technical result of this paper
\begin{lema}\label{keylemma}
Let $\Omega \subset \R^N$ be a bounded smooth domain and $\epsilon \in (0,1)$. For $\beta \in (0,1)$, consider the function
\begin{equation}\label{psi}
\psi(x) = \psi_\beta(x) := (\epsilon + d(x))^\beta \mathbf{1}_{\bar{\Omega}} (x),
\end{equation}
where $d = d_\Omega$ is the distance function to $\partial \Omega$. Then, there exists $\bar{\delta} \in (0, \delta_0)$,
$\beta_0 \in (0, \min \{1, 2\sigma\})$ and a constant $c^* > 0$, depending only on $\Omega, N$ and $\sigma$, such that, for all $\beta \leq \beta_0$
we have
\begin{equation}\label{eqkeylemma}
-\I_\epsilon[\psi](x) \geq c^* (\epsilon + d(x))^{\beta - 2\sigma}, \quad \mbox{for all} \ x \in \bar{\Omega}_{\bar{\delta}}, 
\ \epsilon \in (0, \bar{\delta}).
\end{equation}
\end{lema}

\noindent
{\bf \textit{Proof:}} We start considering $\bar{\delta} < r_0'$ and $x \in \Omega_{\bar{\delta}}$. We split the integral 
\begin{equation*}
\I_\epsilon[\psi](x) = I_0(x) + I_1(x) + I_2(x) + I_3(x), 
\end{equation*}
where
\begin{equation*}
\begin{split}
I_0(x) := & \ \int_{B_{r_0}^c} [\psi(x + z) -\psi(x)]K_\epsilon(z)dz, \\
I_1(x) := & \ \int_{B_{d(x)/2}} [\psi(x + z) - \psi(x)]K_\epsilon(z)dz, \\
I_2(x) := & \ - (\epsilon + d(x))^{\beta - 2\sigma} \int_{(\Omega^c - x)\cap B_{r_0}} K_\epsilon(z)dz \quad\mbox{and} \\
I_3(x) := & \ \int_{(\Omega - x) \cap B_{r_0} \setminus B_{d(x)/2}} [(\epsilon + d(x + z))^\beta - (\epsilon + d(x))^\beta]K_\epsilon(z)dz.
\end{split}
\end{equation*}
In what follows we estimate each $I_i(x), i = 0,1,2,3$.
Since $\psi$ is bounded in $\R^N$ independent of $\epsilon, \beta$ when $\epsilon, \beta \in (0,1)$, we have
\begin{equation}\label{I0}
I_0(x) \leq c_1 r_0^{-2\sigma},
\end{equation}
where $c_1 > 0$ depends only on $\Omega$ and $N$.
For $I_1(x)$, by the symmetry of $K_\epsilon$ we have
\begin{equation*}
I_1(x) = \frac{1}{2} \int_{B_{d(x)/2}} [\psi(x + z) + \psi(x - z) - 2\psi(x)]K_\epsilon(z)dz.
\end{equation*}
Then we consider the
 function $\theta(z)= \psi(x + z) + \psi(x - z) - 2\psi(x)$, which is smooth in $\bar{B}_{d(x)/2}$ and therefore, 
 we can write by Taylor expansion
 \begin{equation*}
\begin{split}
\theta(z) = & \frac{\beta}{2} \Big{[} (\epsilon + d(x + \tilde{z}))^{\beta - 1} \langle D^2d(x + \tilde{z})z, z \rangle \\
& \quad + (\epsilon + d(x - \bar{z}))^{\beta - 1} \langle D^2d(x - \bar{z})z, z \rangle \\
& \quad + (\beta - 1) (\epsilon + d(x + \tilde{z}))^{\beta - 2} |\langle Dd(x + \tilde{z}), z \rangle|^2 \\
& \quad + (\beta - 1) (\epsilon + d(x - \bar{z}))^{\beta - 2} |\langle Dd(x - \bar{z}), z \rangle|^2 \Big{]},
\end{split}
\end{equation*}
where $\tilde{z}, \bar{z} \in B_{d(x)/2}$. With this, since we assume $\beta < 1$, 
by the smoothness of the distance function $d$ inherited by the smoothness of $\partial \Omega$ we have
\begin{equation*}
\theta(z) \leq c_2  (\epsilon + d(x))^{\beta - 1} |z|^2, \quad \mbox{for all} \ z \in B_{d(x)/2},
\end{equation*}
where $c_2=C_\Omega \beta > 0$ depends on the domain, but not on $\epsilon$ or $d(x)$. From this, we get
\begin{equation*}
I_1(x) \leq c_2 (\epsilon + d(x))^{\beta - 1} \int_{B_{d(x)/2}} |z|^2 K_\epsilon(z)dz,
\end{equation*}
and since $K_\epsilon(z) \leq K_0(z)$, we conclude that, for a constant $c_3>0$, we have
\begin{equation}\label{I1}
I_1(x) \leq c_3 \beta (\epsilon + d(x))^{\beta - 2\sigma + 1}.
\end{equation}

Now we address the estimates of $I_2(x)$ and $I_3(x)$. 
For $I_2(x)$, recalling  the change of variables $\Phi_x$, we have 
\begin{equation*}
\Phi_x(B_{r_0'} \times (-r_0', -d(x))) \subset (\Omega^c- x) \cap B_{r_0}.
\end{equation*}
With this, using the change of variables $\Phi_x$ and applying~\eqref{tildeK}, we have
\begin{equation*}
I_2(x)\leq -C_\Omega^{-1} (\epsilon + d(x))^\beta \int_{B_{r_0'} \times (-r_0', -d(x))} K_\epsilon(\xi', s)d\xi' ds.
\end{equation*}
But there exists a constant $c_4 > 0$, depending only on $N$ and $\sigma$, such that
\begin{equation*}
\epsilon^{N + 2\sigma} + |s|^{N + 2\sigma} \leq c_4(\epsilon^{1+ 2\sigma} + |s|^{1 + 2\sigma})^{(N + 2\sigma)/(1 + 2\sigma)},
\end{equation*}
and with this, defining $\rho(\epsilon, s) = (\epsilon^{1 + 2\sigma} + |s|^{1 + 2\sigma})^{1/(1 + 2\sigma)}$, we can write
\begin{equation*}
\begin{split}
I_2(x)\leq & \ -c_5(\epsilon + d(x))^\beta \int \limits_{-r_0'}^{-d(x)} \frac{ds}{\rho(\epsilon, s)^{N + 2\sigma}} 
\int \limits_{B_{r_0'}} \frac{d\xi'}{1 + |\xi'/\rho(\epsilon, s)|^{N + 2\sigma}} \\ 
= & \ -c_5 (\epsilon + d(x))^\beta \int \limits_{-r_0'}^{-d(x)} \frac{ds}{\rho(\epsilon, s)^{1 + 2\sigma}} 
\int \limits_{B_{r_0'}/\rho(\epsilon, s)} \frac{dy}{1 + |y|^{N + 2\sigma}} \\
\leq & \ -c_6 (\epsilon + d(x))^\beta \int \limits_{-r_0'}^{-d(x)} \frac{ds}{\epsilon^{1 + 2\sigma} + |s|^{1 + 2\sigma}}.
\end{split}
\end{equation*}

Finally, making the change $t = -s/(\epsilon + d(x))$, we conclude
\begin{equation}\label{tecnincaI2}
I_2(x)  \leq - c_6 (\epsilon + d(x))^{\beta - 2\sigma} \int \limits_{1 - \tau}^{r_0'/(\epsilon + d(x))}
\frac{dt}{\tau^{1 + 2\sigma} + |t|^{1 + 2\sigma}},
\end{equation}
where $\tau = \epsilon/(\epsilon + d(x)) \in (0,1)$. At this point, taking $\bar{\delta}$ small in order to have 
$\epsilon + d(x) < 2/r_0'$, we find that the interval
$$
(1 - \tau, r_0'/(\epsilon + d(x)))
$$ 
has at least lenght $1$. 
Hence, we conclude the existence of $c_7 > 0$, depending only on $\Omega, N$ and $\sigma$, such that
\begin{equation}\label{I2}
I_2(x)  \leq - c_7 (\epsilon + d(x))^{\beta - 2\sigma}. 
\end{equation}

It remains to estimate $I_3(x)$. Defining $D_+(x) = \{ z : d(x + z) \geq d(x)\}$, we clearly have
\begin{equation*}
I_3(x)  \leq \int \limits_{(\Omega - x) \cap D_+(x) \cap B_{r_0} \setminus B_{d(x)/2} } [(\epsilon + d(x + z))^\beta - (\epsilon + d(x))^\beta]
K_\epsilon(z)dz,
\end{equation*}
and since 
\begin{equation*}
(\Omega - x) \cap B_{r_0} \cap D_+ \subset \Phi_x(B_{r_0} \times (0, r_0)),
\end{equation*}
we have 
\begin{equation*}
\begin{split}
I_3(x) \leq & \ \int \limits_{\Phi_x(B_{r_0} \times (0, r_0)) \setminus B_{d(x)/2}} [(\epsilon + d(x + z))^\beta - (\epsilon + d(x))^\beta]
 K_\epsilon(z)dz \\ 
= & \ C_\Omega (\epsilon + d(x))^\beta \int \limits_{\Phi_x(B_{r_0} \times (0, r_0)) \setminus B_{d(x)/2}} 
\Big{[} \Big{(} \frac{\epsilon + d(x + z)}{\epsilon + d(x)} \Big{)}^\beta - 1 \Big{]} K_\epsilon(z)dz,
\end{split}
\end{equation*}

Thus, making a change of variables we have
\begin{equation*}
I_3(x) \leq C_\Omega (\epsilon + d(x))^\beta \int \limits_{B_{r_0} \times (0, r_0) \setminus \Phi_x^{-1}(B_{d(x)/2})}  
\Big{[} \Big{(} \frac{\epsilon + d(x + \Phi_x(\xi))}{\epsilon + d(x)} \Big{)}^\beta - 1 \Big{]} \tilde{K}_\epsilon(\xi)d\xi.
\end{equation*}

Since $\Phi_x$ is a diffeomorphism, there exists a constant $c_8 > 0$ such that
\begin{equation*}
d(x + \Phi_x(\xi', s)) \leq d(x + \Phi_x(0, s)) + c_8 |\xi'| = d(x) + s + c_8 |\xi'|
\end{equation*}
and a constant $\lambda \in (0,1)$ small, depending only on the smoothness of $\partial \Omega$, such 
that $B_{\lambda d(x)} \subset \Phi_x^{-1}(B_{d(x)/2})$. Using this and~\eqref{tildeK} we arrive to
\begin{equation*}
\begin{split}
I_3(x) \leq & \ C_\Omega (\epsilon + d(x))^\beta \int \limits_{B_{r_0} \times (0, r_0) \setminus B_{\lambda d(x)}}  
[(1 + c_8|\xi/(\epsilon + d(x))|)^\beta - 1] K_\epsilon(\xi)d\xi \\
= & \ C_\Omega (\epsilon + d(x))^{\beta - 2\sigma} \int \limits_{(\epsilon + d(x))^{-1} B_{r_0} \setminus B_{\lambda d(x)}}
[(1 + c_8|y|)^\beta - 1] K_\tau(y)dy \\
\leq & \ C_\Omega (\epsilon + d(x))^{\beta - 2\sigma} \int_{\lambda (1 - \tau)}^{+\infty} \frac{[(1 + c_8 t)^\beta - 1] t^{N - 1} dt}{\tau^{1 + 2\sigma} + t^{N + 2\sigma}}.
\end{split}
\end{equation*}

At this point, we remark that for each $M > 2$, we have
\begin{equation*}
\int_{M}^{+\infty} \frac{[(1 + c_8 t)^\beta - 1] t^{N - 1} dt}{\tau^{1 + 2\sigma} + t^{N + 2\sigma}} \leq c_9 M^{\beta - 2\sigma},
\end{equation*}
where $c_9 > 0$ depends only on $N, \sigma$ and $\Omega$. On the other hand, for each $M > 2$ there exists $\beta = \beta(M) > 0$  small such that
\begin{equation*}
\int_{\lambda (1 - \tau)}^{M} \frac{[(1 + c_8 t)^\beta - 1] t^{N - 1} dt}{\tau^{1 + 2\sigma} + t^{N + 2\sigma}} 
\leq C_\Omega^{-1} c_7/2,
\end{equation*}
where $c_7 > 0$ is the constant arising in~\eqref{I2}. From the last two estimates, we conclude that for each $M > 2$, there exists $\beta$
small such that
\begin{equation}\label{I3}
I_3(x) \leq  c_7(\epsilon + d(x))^{\beta - 2\sigma}/2 + c_{10}M^{\beta - 2\sigma},
\end{equation}
where $c_{10} > 0$ depends only on $N, \sigma$ and $\Omega$.
Putting together ~\eqref{I0},~\eqref{I1},~\eqref{I2} and~\eqref{I3}, and fixing $M = \max \{ 2, r_0\}$, we have
\begin{equation*}
\I_\epsilon[\psi](x)  \leq (\epsilon + d(x))^{\beta - 2\sigma} (-c_7 /2 + c_2 \beta (\epsilon + d(x))) + c_{11} r_0^{-2\sigma},
\end{equation*}
where $c_{11} > 0$ depends only on $N, \sigma$ and $\Omega$. Hence, fixing $\beta > 0$ smaller if it is necessary, we can write
\begin{equation*}
\I_\epsilon[\psi](x)  \leq -c_7 (\epsilon + d(x))^{\beta - 2\sigma}/4 + c_{11} r_0^{-2\sigma}. 
\end{equation*}

Finally, taking $\epsilon + d(x)$ small in terms of $c_7, c_{11}, r_0, \beta$ and $\sigma$ (and therefore, depending only on $N, \sigma$ and $\Omega$),
we conclude~\eqref{eqkeylemma}, where $c^* = c_7/8$.
\qed


\medskip

The last lemma allows us to provide the following control of the discontinuity at the boundary.
\begin{prop}\label{propaprioriborde}
Let $\epsilon \in (0,1)$ and $u_\epsilon$ the solution of~\eqref{eq}-\eqref{exteriordata}. 
Let $\bar{\delta} > 0$ and $\beta_0 \in (0, \min \{ 1, 2\sigma\})$
as in Lemma~\ref{keylemma}. Then, for each $d_0 \in (0,\bar{\delta})$, there exists $C_0 > 0$ satisfying
\begin{eqnarray*}
|u_\epsilon(x)| \leq C_0(\epsilon + d(x))^{\beta_0} \quad \mbox{for all} \ x \in \bar{\Omega}_{d_0}.
\end{eqnarray*}
The constant $C_0$ depends on $\beta_0, d_0$, $\sigma$ and $\Omega$.
\end{prop}

\noindent
{\bf \textit{Proof:}} Let $\beta_0$ as in Lemma~\ref{keylemma}, $\psi = \psi_{\beta_0}$ as in~\eqref{psi} 
and consider the function
\begin{equation*}
\zeta(x) = \min \{ \psi(x), (\epsilon + d_0)^{\beta_0} \}, \quad x \in \R^N.
\end{equation*}
Observing that $\zeta = \psi$ in $\Omega_{d_0} \cup \Omega^c$ and $\psi \geq \zeta$ in $\R^N$, we easily conclude that
\begin{equation*}
\I_\epsilon[\zeta](x) \leq \I_\epsilon[\tilde{\zeta}](x), \quad \mbox{for all} \ x \in \bar{\Omega}_{d_0},
\end{equation*}
and using Lemma~\ref{keylemma} we get
\begin{equation*}
-\I_\epsilon[\zeta](x) \geq c^* (\epsilon + d(x))^{\beta - 2\sigma} \quad \mbox{for all} \ x \in \bar{\Omega}_{d_0}.
\end{equation*}

Let ${C} > 0$ be the constant in Proposition~\ref{cotauniformesolucion}
and define the function $\tilde{z}_+ = (C d_0^{-\beta} + 2^\sigma {c^*}^{-1}) ||f||_\infty \zeta$. By  construction of $\tilde{z}_+$, we have
\begin{equation*}
-\I_\epsilon [\tilde{z}_+] \geq ||f||_\infty \quad \mbox{in} \ \bar{\Omega}_{d_0}; 
\qquad \mbox{and} \qquad \tilde{z}_+ \geq u_\epsilon \quad \mbox{in} \ \bar{\Omega}_{d_0}^c,
\end{equation*}
and therefore, applying the comparison principle, we conclude $u_\epsilon \leq \tilde{z}_+$ in $\bar{\Omega}_{d_0}$.
Similarly, we can conclude the function $\tilde{z}_- = -\tilde{z}_+$ satisfies $\tilde{z}_- \leq u_\epsilon$ in $\bar{\Omega}_{d_0}$,
from which we get the result.
%
\qed


\section{Boundary Equicontinuity.}
\label{seccionND}

In this section we establish the boundary equicontinuity of the family of solutions $\{ u_\epsilon \}_{\epsilon \in (0,1)}$ of
problem~\eqref{eq}-\eqref{exteriordata}.
The main result of this section  is the following
\begin{teo}\label{teoborde}
Let $\epsilon \in (0,1)$ and $u_\epsilon$ be the solution to~\eqref{eq}-\eqref{exteriordata}. There exists a modulus of 
continuity $m_0$ depending only on $N, \sigma, f$ and $\Omega$, such that 
\begin{equation*}
|u_\epsilon(x) - u_\epsilon(y)| \leq m_0(|x - y|) \quad \mbox{for all} \ x,y \in \bar{\Omega}_{\bar{\delta}},
\end{equation*}
with $\bar{\delta} > 0$ given in Lemma~\ref{keylemma}.
\end{teo}

The idea of the proof is based on the fact $w(x) = u(x+y)-u(x)$, where $y$ is fixed, satisfies an equation (near the boundary) 
for which the comparison principle holds. Using this, we get the result
constructing a barrier to this problem, independent of $\epsilon$ and associated to $m$ in Theorem~\ref{teoborde}.

In what follows we discuss the  precise elements on the proof. We consider $y \in \R^N$ with $0 < |y| < \bar{\delta}/2$, 
with $\bar{\delta}$ as in Lemma~\ref{keylemma}. Define the sets
\begin{equation*}\label{defOyU}
\mathcal{O} = \mathcal{O}(y) := \Omega \setminus \bar{\Omega}_{|y|}, \quad \mathcal{U} = \mathcal{U}(y) := \{x \in \R^N : -|y| \leq d_\Omega(x) < |y| \}.
\end{equation*}
and the function
\begin{equation}\label{defw}
w(x) = w_{y, \epsilon}(x) := u_\epsilon(x + y) - u_\epsilon(x), \ x \in \R^N. 
\end{equation}

Notice that $w \equiv 0$ in $\R^N \setminus (\bar{\mathcal{O}} \cup \U)$ and, by Proposition~\ref{propaprioriborde}, there exists $C_0, \beta_0 > 0$ 
such that $|w(x)| \leq C_0(\epsilon + |y|)^{\beta_0}$ for all $x \in \U$.
Since we have that $w$ satisfies
\begin{eqnarray*}
- \I_\epsilon[w] (x) = f(x + y) - f(x) \quad \mbox{for all } \ x \in \bar{\mathcal{O}},
\end{eqnarray*}
denoting by $m_f$  the modulus of continuity of $f$, we conclude that $w \in C(\bar{\mathcal{O}})$ satisfies the inequality
\begin{equation}\label{eqw}
- \I_\epsilon[w] (x) \leq m_f(|y|) \quad \mbox{in} \ \bar{\mathcal{O}},
\end{equation}
and the exterior inequality
\begin{equation}\label{exteriordataw}
w(x) \leq C_0 (\epsilon + |y|)^{\beta_0} \mathbf{1}_{\U}(x) \quad \mbox{in} \ \bar{\mathcal{O}}^c.
\end{equation}

Let $\zeta$ and $\eta$ the functions defined as
\begin{equation*}
\begin{split}
\zeta(x) = & \min \{ (\epsilon + \bar{\delta} - |y|)^\epsilon, (\epsilon + d_\Omega(x) - |y|)^\epsilon \} \mathbf{1}_{\bar{\mathcal{O}}}(x)\quad\mbox{and} \\
\eta(x) = & C_0(\epsilon + |y|)^{\beta_0} \mathbf{1}_{\U}(x),
\end{split}
\end{equation*}
and consider the function
\begin{eqnarray}\label{defW}
W(x) = \eta(x) + A m(|y|) \zeta(x),
\end{eqnarray}
where $A > 0$ and $m$ is a modulus of continuity satisfying $m(|y|) \geq m_f(|y|)$. We have the following
\begin{prop}\label{propWsupersol}
There exists $A > 0$ large, depending on $\Omega, N$ and $\sigma$, such that
\begin{equation*}
-\I_\epsilon[W](x) \geq m_f(|y|), \quad \mbox{for all} \ x \in \bar{\mathcal{O}},
\end{equation*}
for all $\epsilon \in (0,\bar{\delta})$, with $\bar{\delta}$ given in Lemma~\ref{keylemma}.
\end{prop}

\noindent
{\bf \textit{Proof:}} Without loss of 
generality we may assume the existence of a number $0 < \alpha < \min \{1, \beta_0\}$ and a constant $c_1$ such that
\begin{eqnarray}\label{m>potencia}
m(t) \geq c_1 t^{\alpha}, \quad \mbox{for all} \ t \geq 0.
\end{eqnarray}
By linearity of $\I_\epsilon$, we have 
\begin{equation*}
\I_\epsilon[W](x) = \I_\epsilon[\eta](x) + Am(|y|) \I_\epsilon[\zeta](x).
\end{equation*}
Thus, we may estimate each term in the right-hand side separately. 

\medskip
\noindent
\textsl{1.- Estimate for $\I_\epsilon[\zeta](x)$:}
We first notice that for $x \in \Omega$ with $|y| \leq d_\Omega(x) \leq \bar{\delta}$ we can write
\begin{equation*}
\zeta(x) = (\epsilon + d_\Omega(x) - |y|)^\epsilon \mathbf{1}_{\bar{\mathcal{O}}} 
= (\epsilon + d_{\mathcal{O}}(x))^\epsilon \mathbf{1}_{\bar{\mathcal{O}}}(x).
\end{equation*}
Then, applying Lemma~\ref{keylemma}, for all $\epsilon$ small we have
\begin{equation*}
-\I_\epsilon[\zeta](x) \geq c^* (\epsilon + d(x) - |y|)^{\epsilon - 2\sigma}, \quad 
\mbox{for all} \ x \in \bar{\Omega}_{\bar{\delta}} \cap \bar{\mathcal{O}}, 
\end{equation*}
for some $c^* > 0$ not depending on $d(x), |y|$ or $\epsilon$. In fact, for all $\epsilon \in (0,1)$ the term $(\epsilon + d(x) - |y|)^{-\epsilon}$
is bounded below by a strictly positive constant, independent of $\epsilon$, driving us to
\begin{equation}\label{Izetaborde1}
-\I_\epsilon[\zeta](x) \geq c^* (\epsilon + d(x) - |y|)^{- 2\sigma}, \quad 
\mbox{for all} \ x \in \bar{\Omega}_{\bar{\delta}} \cap \bar{\mathcal{O}}.
\end{equation}

On the other hand, when $x \in \Omega \setminus \bar{\Omega}_{\bar{\delta}}$, for all $\epsilon \in (0,1)$ we have
\begin{equation*}
\I_\epsilon[\zeta](x) \leq -(\epsilon + \bar{\delta} - |y|)^\epsilon \int \limits_{(\Omega \setminus \Omega_{|y|})^c - x} K_\epsilon(z)dz
\leq -\epsilon^\epsilon \int \limits_{\Omega^c - x} K_1(z)dz,
\end{equation*}
and therefore, there exists $c_2 > 0$, not depending on $\epsilon, d(x)$ or $|y|$, such that
\begin{equation*}\label{Izetainterior}
\I_\epsilon[\zeta](x) \leq -c_2, \quad \mbox{for all} \ x \in \Omega \setminus \bar{\Omega}_{\bar{\delta}}.
\end{equation*}
Since $|y| \leq \bar{\delta}/2$, making 
$c^*$ smaller if  necessary, the last inequality and~\eqref{Izetaborde1} drives us to
\begin{equation}\label{Izetaborde}
-\I_\epsilon[\zeta](x) \geq c^* (\epsilon + d(x) - |y|)^{- 2\sigma}, \quad 
\mbox{for all} \ x \in \bar{\mathcal{O}},  
\end{equation}

\medskip
\noindent
\textsl{2.- Estimate for $\I_\epsilon[\eta](x)$:} By its very definition, for $x \in \bar{\mathcal{O}}$ we have
\begin{equation}\label{Ietadef}
\I_\epsilon[\eta](x) = C_0 (\epsilon + |y|)^{\beta_0} \int \limits_{\mathcal{U} - x} K_\epsilon(z)dz.
\end{equation}
We start considering the case $x \in \Omega \setminus \Omega_{\bar{\delta}}$, where  we have $\mathrm{dist}(x, \U) \geq \bar{\delta}/2$ and then, there exists a constant $c_3>0$ depending only on $\bar{\delta}$
(which in turn depends only on the smoothness of the domain), such that $K_\epsilon(z) \mathbf{1}_{\U - x} \leq c_3$. 
Using this, we have
\begin{equation*}
\I_\epsilon[\eta](x) \leq c_3 (\epsilon + |y|)^{\beta_0} \int \limits_{\U - x} dz.
\end{equation*}

By the boundedness of $\Omega$, there exists $c_4 > 0$ depending only on $N$ such that $\mathrm{Vol}(\U - x) \leq c_4 |y|$.
Using this and~\eqref{m>potencia}, we conclude that
\begin{equation}\label{Ietacaso1}
\I_\epsilon[\eta](x) \leq c_5 m(|y|),
\end{equation}
where $c_5 > 0$ depends only on $N, \sigma$ and $\Omega$. 

Now we deal with the case $x \in \mathcal{O} \cap \Omega_{\bar{\delta}}$ (notice that in this case we are assuming $d_\Omega(x) > |y|$).
Using~\eqref{Ietadef} and recalling the change of variables $\Phi_x$ introduced in~\eqref{Phi}, we can write
\begin{equation*}
\I_\epsilon[\eta](x) \leq C_0 (\epsilon + |y|)^{\beta_0} \Big{(} \int \limits_{(\U - x) \setminus B_{r_0}} K_\epsilon(z)dz 
+ \int \limits_{(\U - x) \cap \mathcal{R}_x} K_\epsilon(z)dz \Big{)},
\end{equation*}
where $\mathcal{R}_x$ was defined at the beginning of Section \S 3.
Using a similar analysis as the one leading to~\eqref{Ietacaso1}, there exists a universal constant $c_6 > 0$ such that
\begin{equation*}
\I_\epsilon[\eta](x) \leq  c_6 (\epsilon + |y|)^{\beta_0} \Big{(} |y| 
+ \int \limits_{(\U - x) \cap \mathcal{R}_x} K_\epsilon(z)dz\Big{)}.
\end{equation*}

Now, we have that $(\U - x) \cap \mathcal{R}_x = \Phi_x(B_{R_0} \times (-d(x) - |y|, |y| - d(x)))$, and therefore, applying the change of variables
$\Phi_x$ and the estimate~\eqref{tildeK}, we arrive to
\begin{equation*}
\I_\epsilon[\eta](x) \leq  c_7 (\epsilon + |y|)^{\beta_0} \Big{(} |y| 
+ \int \limits_{B_{R_0} \times (-d(x) - |y|, |y| - d(x))} \frac{d\xi' ds}{\epsilon^{N + 2\sigma} + |(\xi', s)|^{N + 2\sigma}} \Big{)},
\end{equation*}
and from this, using a similar argument as the one leading to~\eqref{tecnincaI2} to treat the last integral term, and applying~\eqref{m>potencia}, 
we conclude that
\begin{equation}\label{Ieta0}
\I_\epsilon[\eta](x) \leq  c_8 (\epsilon + |y|)^{\beta_0} \Big{(} m(|y|) 
+ \int \limits_{d(x) - |y|}^{d(x) + |y|} \frac{ds}{\epsilon^{1 + 2\sigma} + |s|^{1 + 2\sigma}} \Big{)}.
\end{equation}

Now, the core of this estimate is the computation of the last integral. Denoting
\begin{equation*}
I(x) := (\epsilon + |y|)^{\beta_0} 
\int \limits_{d(x) - |y|}^{d(x) + |y|} \frac{ds}{\epsilon^{1 + 2\sigma} + |s|^{1 + 2\sigma}}, 
\end{equation*}
we claim the existence of a constant $c_9 > 0$ not depending on $\epsilon, d(x)$ or $|y|$ such that
\begin{equation}\label{I(x)}
I(x)  \leq  c_9 m(|y|) (\epsilon + d(x) - |y|)^{- 2\sigma}.
\end{equation}
We get this estimate considering various  cases. When $|y| \leq \epsilon$ and $d(x) - |y| \leq 2\epsilon$ we write
\begin{equation*}
I(x) = (\epsilon + |y|)^{\beta_0} \epsilon^{-2\sigma} \int \limits_{(d(x) - |y|)/\epsilon}^{(d(x) + |y|)/\epsilon}K_1(z)dz
\leq 2^{\beta_0 + 1} \epsilon^{\beta_0-2\sigma - 1} |y|, 
\end{equation*}
and using that $m(|y|) \geq |y|^\alpha$ for some $\alpha \in (0, \beta_0)$, we have
\begin{equation*}
\begin{split}
I(x) \leq & \ 2^{\beta_0 + 1} m(|y|) \epsilon^{\beta_0 - 1} \epsilon^{-2\sigma} |y|^{1 - \alpha} \\
\leq & \  2^{\beta_0 + 1} 3^{2\sigma} m(|y|) \epsilon^{\beta_0 - \alpha}(\epsilon + d(x) - |y|)^{- 2\sigma},
\end{split}
\end{equation*}
and from this, we conclude 
\begin{equation}\label{Ixcaso1}
I(x) \leq c_{10} \epsilon^{\beta_0 - \alpha} m(|y|) (\epsilon + d(x) - |y|)^{- 2\sigma},
\end{equation}
for some constant $c_{10} > 0$.

When $|y| \leq \epsilon$ and $d(x) - |y| > 2\epsilon$, we have
\begin{equation*}
I(x) \leq 2^{\beta_0} \epsilon^{\beta_0} \int_{d(x) - |y|}^{d(x) + |y|} |z|^{-(1 + 2\sigma)}dz 
\leq 2^{\beta_0} \epsilon^{\beta_0} (d(x) - |y|)^{-(1 + 2\sigma)} |y|,
\end{equation*}
and using that $m(|y|) \geq |y|^\alpha$, we arrive to
\begin{equation*}
I(x) \leq 2^{\beta_0 - 1} m(|y|) \epsilon^{\beta_0 - \alpha} (d(x) - |y|)^{-2\sigma} 
\leq 2^{\beta_0 - 1 + 2\sigma}  m(|y|) \epsilon^{\beta_0 - \alpha} (\epsilon + d(x) - |y|)^{-2\sigma},
\end{equation*}
concluding the same estimate~\eqref{Ixcaso1}.

In the case $|y| > \epsilon$ and $d(x) - |y| \leq 2\epsilon$, performing the change $\xi = z/\epsilon$ in the integral defining 
$I(x)$, we have
\begin{equation*}
I(x) \leq (\epsilon + |y|)^{\beta_0} \epsilon^{-2\sigma} ||K_1||_{L^1} 
\leq ||K_1||_{L^1} 3^{2\sigma} 2^{\beta_0} |y|^{\beta_0 - \alpha} m(|y|) (\epsilon + d(x) - |y|)^{-2\sigma},
\end{equation*}
and therefore we conclude
\begin{equation}\label{Ixcaso3}
I(x) \leq C |y|^{\beta_0 - \alpha} m(|y|) (\epsilon + d(x) - |y|)^{-2\sigma}. 
\end{equation}

Finally, in the case $|y| > \epsilon$ and $d(x) - |y| > 2\epsilon$ we have
\begin{equation*}
I(x) \leq (\epsilon + |y|)^{\beta_0} \epsilon^{-2\sigma} \int \limits_{(d(x) - |y|)/\epsilon}^{(d(x) + |y|)/\epsilon} K_1(z)dz
\leq 2^{\beta_0 -1} \sigma^{-1} |y|^{\beta_0} (d(x) - |y|)^{-2\sigma}, 
\end{equation*}
from which we arrive to~\eqref{Ixcaso3}. From~\eqref{Ixcaso1} and~\eqref{Ixcaso3} we arrive to~\eqref{I(x)}. 
Hence, there exists $c_{11} > 0$ depending only on $N, \Omega$ and $\sigma$ such that
\begin{equation*}
-\I_\epsilon[\eta](x) \geq - c_{11} m(|y|) \Big{(} (\epsilon + |y|)^{\beta_0} +  (\epsilon + d(x) - |y|)^{- 2\sigma} \Big{)},
\end{equation*}
for $x \in \mathcal{O} \cap \Omega_{\bar{\delta}}$. Taking this inequality and~\eqref{Ietacaso1}, since $|y| \leq \bar{\delta}/2$
there exists a constant $c_{12} > 0$ such that
\begin{equation}\label{Ietaborde}
-\I_\epsilon[\eta](x) \geq - c_{12}m(|y|) \Big{(} (\epsilon + |y|)^{\beta_0} +  (\epsilon + d(x) - |y|)^{- 2\sigma} \Big{)}, 
\end{equation}
for all $x \in \bar{\mathcal{O}}$, where the estimate for $x \in \partial \mathcal{O}$ is valid by Lemma~\ref{lemahastaelborde}.

\medskip
\noindent
\textsl{3.- Conclusion:} For each $x \in \bar{\mathcal{O}}$, by~\eqref{Izetaborde} and~\eqref{Ietaborde} we have
\begin{equation*}
-\I_\epsilon[W](x) \geq \Big{[} ( A c^* - c_{12} )  (\epsilon + d(x) - |y|)^{-2\sigma} - c_{12}(\epsilon + |y|)^{\beta_0} \Big{]}m(|y|),
\end{equation*}
and therefore, by taking $A$ large in terms of $N,  \sigma, c_{12}, c^*$ and $\mathrm{diam}(\Omega)$, we conclude by the choice of $m$ that
\begin{equation*}
-\I_\epsilon[W](x) \geq m(|y|) \geq m_f(|y|), \quad \mbox{for all} \ x \in \bar{\mathcal{O}},
\end{equation*}
and the proof follows.
\qed

\medskip

This proposition allows us to give the

\noindent
{\bf \textit{Proof of Theorem~\ref{teoborde}:}} 
Since $w$ defined in~\eqref{defw} satisfies problem~\eqref{eqw}-\eqref{exteriordataw} and recalling $W$ defined in~\eqref{defW}, 
by Proposition~\ref{propWsupersol} and the form of $W$ in $\bar{\mathcal{O}}^c$, we can use the comparison principle to 
conclude that $w \leq W$ in $\bar{\mathcal{O}}$. This means that
\begin{equation*}
u_\epsilon(x + y) - u_\epsilon(x) = w(x) \leq W(x) \leq c_1 Am(y), \quad x \in \bar{\Omega}_y,
\end{equation*}
for some constant $c_1 > 0$. Since a similar lower bound can be stated, by the arbitrariness of $y$ 
we conclude the result with $m_0 = c_1Am$.
\qed


\section{Proof of Theorem~\ref{teo1}.}
\label{pruebamodulointerior}

Consider $\bar{\delta}$ as in Lemma~\ref{keylemma}, let $y \in \R^N$ such that $|y| \leq \bar{\delta}/8$ and consider the sets
\begin{eqnarray*}
& \Sigma_1  = \overline{(\Omega - y) \cup \Omega}, & \ \Sigma_2 = \Omega \cap (\Omega - y), \\
& \Sigma_3 = \Sigma_1 \setminus \Sigma_2 \quad\mbox{and}\quad
& \ \Sigma_4 = (\Omega \setminus \bar{\Omega}_{\bar{\delta}/2}) \cup ((\Omega \setminus \bar{\Omega}_{\bar{\delta}/2}) - y). 
\end{eqnarray*}

Notice that $\Sigma_4 \subset \Sigma_2 \subset \Sigma_1$. In addition, notice that if $z \in \Sigma_3$, then $z + y$ and $z$ 
cannot be simultaneously in $\Omega$. We also have 
$$
|\mathrm{dist}(z, \partial \Omega)|, |\mathrm{dist}(z + y, \partial \Omega)| \leq |y|
$$ 
for each $z \in \Sigma_3$. Finally, observe that if $x \in \Sigma_2 \setminus \Sigma_4$, then $x, x + y \in \Omega_{\bar{\delta}}$.
Thus, considering $w$ as in~\eqref{defw}, by Proposition~\ref{propaprioriborde} we can assure the existence of $C_0, \beta_0 > 0$
such that
\begin{equation*}
w \leq C_0 (\epsilon + |y|)^{\beta_0} \quad \mbox{in} \ \Sigma_3,
\end{equation*}
and by Theorem~\ref{teoborde} we have 
\begin{equation*}
w \leq m_0(|y|) \quad \mbox{in} \ \Sigma_2 \setminus \bar{\Sigma}_4.
\end{equation*}
Now, consider the funtion
\begin{equation*}
Z(x) = A m_0(|y|) \mathbf{1}_{\Sigma_2}(x) + C_0(\epsilon + |y|)^{\beta_0} \mathbf{1}_{\Sigma_3}(x),
\end{equation*}
where $A > 0$ is a constant to be fixed later. 
Notice that for each $x \in \bar{\Sigma}_3$, we have
\begin{equation}\label{IZ}
\I_\epsilon[Z](x) = C_0 (\epsilon + |y|)^{\beta_0} \int \limits_{\Sigma_3 - x} K_\epsilon(z)dz 
- Am_0(|y|) \int \limits_{\Sigma_2^c - z} K_\epsilon(z)dz.
\end{equation}

At this point, we remark that there exists a constant $c_1 > 0$, independent of $\epsilon, y$ or $x$, such that 
\begin{equation*}
\int \limits_{\Sigma_2^c - z} K_\epsilon(z)dz \geq c_1.
\end{equation*}
On the other hand, since $\mathrm{dist}(x, \Sigma_3) \geq \bar{\delta}/2$ we have $K_\epsilon(z) \mathbf{1}_{\Sigma_3 - x} \leq c_2$ for some
 constant $c_2 > 0$, and by the boundedness of $\Omega$, $\mathrm{Vol}(\Sigma_3 - x) \leq c_3|y|$ for some $c_3>0$. Using these facts on~\eqref{IZ}
and applying~\eqref{m>potencia}, we arrive to
\begin{equation*}
\I_\epsilon[Z](x) \leq (c_4 (\epsilon + |y|)^{\beta_0} - c_1A )m_0(|y|).
\end{equation*}
Thus, taking $A$ large in terms of $c_1, c_4$, we conclude that $-\I_\epsilon[Z] \geq m_f(|y|)$ in $\bar{\Sigma}_4$.

By the very definition of $w$, we have
\begin{equation*}
- \I_\epsilon[w] = f(x + y) - f(x), \quad \mbox{for} \ x \in \bar{\Sigma}_4. 
\end{equation*}
Then we have that $-\I_\epsilon[Z] \geq -\I_\epsilon[w]$ in $\bar{\Sigma}_4$ and by definition of $W$ and the bounds of $w$ in $\bar{\Sigma}_4^c$
stated above, we conclude that $w \leq W$ in $\bar{\Sigma}_4^c$. Using the comparison principle, we conclude $w \leq W$ in $\bar{\Sigma}_4$.
A similar argument states the inequality $-W \leq w$ and the result follows.
\qed

\section{Further Results.}
\label{remarks}

\subsection{Fully Nonlinear Equations.}
The result obtained in Theorem~\ref{teoborde} can be readily extended to a certain class of fully nonlinear equations. For example, consider 
two sets of indices $\mathcal{A}, \mathcal{B}$ and a two parameter family of radial continuous functions $a_{\alpha \beta}: \R^N \to \R$ satisfying 
the \textsl{uniform ellipticity condition}
\begin{eqnarray}\label{elipticidadnolineal}
\lambda_1 \leq a_{\alpha \beta}(z) \leq \lambda_2, \quad \forall \alpha \in \mathcal{A}, \beta \in \mathcal{B}, z \in \R^N 
\end{eqnarray}
for certain constants $\lambda_1, \lambda_2$ such that $0 < \lambda_1 < \lambda_2 < +\infty$. Let us denote
$$
K_{\alpha \beta, \epsilon}(z) := \frac{a_{\alpha \beta} (z)}{\epsilon^{n + 2\sigma} + |z|^{n + 2\sigma}}
$$
and with this, for a suitable function $u$ and $x \in \R^N$, define the linear operators
$$
L_{\alpha \beta, \epsilon}[u](x) := \int_{\R^N} \delta(u,x,z) K_{\alpha \beta, \epsilon}(z) dz
$$
and the corresponding \textsl{Isaacs Operator}
$$
I_\epsilon[u](x) = \inf_{\alpha \in \mathcal{A}} \sup_{\beta \in \mathcal{B}} L_{\alpha \beta, \epsilon}[u](x).
$$

Under these definitions, we may consider the corresponding nonlinear equation
\begin{eqnarray}\label{eqnonlinear}
\left \{ \begin{array}{rll} - I_\epsilon[u] =& f &\quad \mbox{in} \ \Omega \\ u =& 0 &\quad \mbox{in} \ \Omega^c. \end{array} \right .
\end{eqnarray}
Existence and uniqueness of a pointwise solution $u_\epsilon$ to~\eqref{eqnonlinear}, which is continuous in $\bar{\Omega}$ can be obtained in a very similar 
way as in the linear case, and Proposition~\ref{lemahastaelborde} can be adapted to this nonlinear setting. This allows us to use the 
comparison principle stated in Proposition~\ref{comparison} as well. 

The lack of linearity can be handled with the positive homogeneity of these operators and the so called \textsl{extremal operators}
\begin{eqnarray*}
\mathcal{M}^+_\epsilon[u](x) = \sup_{\alpha \in \mathcal{A}, \beta \in \mathcal{B}} L_{\alpha \beta, \epsilon}[u](x), \quad
\mathcal{M}^-_\epsilon[u](x) = \inf_{\alpha \in \mathcal{A}, \beta \in \mathcal{B}} L_{\alpha \beta, \epsilon}[u](x), 
\end{eqnarray*}
since, for two functions $u_1, u_2$ and $x \in \R^N$, these operators satisfy the fundamental inequality
\begin{eqnarray*}
\mathcal{M}^-_\epsilon[u_1 - u_2](x) \leq \I_\epsilon[u_1](x) - \I_\epsilon[u_2](x) \leq \mathcal{M}^+_\epsilon[u_1 - u_2](x).
\end{eqnarray*}

A priori estimates for the solution as it is stated in 
Proposition~\ref{propaprioriborde} can be found using the same barriers given in the proof of that proposition, 
as the following useful estimates hold:
For each 
$\alpha \in \mathcal{A}$, $\beta \in \mathcal{B}$, $D \subset \R^N$ 
\begin{eqnarray*}
\int_{D} h K_{\alpha \beta, \epsilon}(z)dz \leq \lambda_1 \int_{D} h K_{\epsilon}(z)dz,
\end{eqnarray*}
for all $h: D \to \R$ bounded nonnegative function, and
\begin{eqnarray*}
\int_{D} h K_{\alpha \beta, \epsilon}(z)dz  \leq \lambda_2 \int_{D} h K_{\epsilon}(z)dz
\end{eqnarray*}
for all $h: D \to \R$ bounded nonpositive function. 

Using these inequalities and \eqref{elipticidadnolineal}, we can use the same barriers appearing in the proof
of Theorem~\ref{teo1} (Theorem~\ref{teoborde} included) and get similar result. Moreover, the same modulus of continuity for the linear case
can be obtained in this nonlinear framework, up to a factor depending on $\lambda_1$ and $\lambda_2$.

\subsection{Parabolic Equations.}

Let $T > 0$, $f: \bar{\Omega} \times [0,T] \to \R$ be a continuous function. A  result similar to Theorem~\ref{teoborde} can be readily obtained for 
the parabolic nonlocal equation
\begin{eqnarray}\label{parabolico}
\left \{ \begin{array}{rll} u_t - \I_\epsilon[u] =& f &\quad \mbox{in} \ \Omega \times [0,T), \\ 
u(x,t) =& 0 &\quad \mbox{in} \ \Omega^c \times [0,T), \\
u(x,0) =& 0 &\quad \mbox{in} \ \bar{\Omega}. \end{array} \right . 
\end{eqnarray}

Similar problem is adressed by the authors in~\cite{Felmer-Topp} for the Cauchy problem in all $\R^N$. Inspired by techniques 
used by Ishii in~\cite{Ishii2},
a modulus of continuity in time can be derived once a modulus of continuity in space is found. So, the key fact is the modulus in $x$ 
and this can be obtained in the parabolic setting noting that Theorem~\ref{teoborde} readily applies considering equations with the form
$$
\lambda u - \I_\epsilon[u] = f \quad \mbox{in} \ \Omega
$$
for $\lambda > 0$ and that each time $Z(x)$ is a suitable barrier for this problem, then the function $(x,t) \mapsto e^t Z(x)$ plays the role of
a barrier for the evolution problem~\eqref{parabolico}.

\subsection{Convergence Issues.}
The proof of Corollary~\ref{corolarioteo} is  standard in the viscosity sense, once the uniform convergence is stated.
However, following the ideas of Cort\'azar, Elgueta and Rossi in ~\cite{C-E-Rossi}, and also in ~\cite{Felmer-Topp}, under stronger
assumptions over the regularity of $u$ in Corollary~\ref{corolarioteo}, we can find a rate of convergence.

\begin{teo}\label{tasa1}
Let $f$, $u_\epsilon$  and $u$ as in Corollary~\ref{corolarioteo} and assume $u \in C^{2\sigma + \gamma}(\bar{\Omega})$ for some $\gamma > 0$. Then,
$$
||u_\epsilon - u||_{L^{\infty}(\bar{\Omega})} \leq C\epsilon^{\gamma_0}
$$ 
for some $0 < \gamma_0 \leq \min \{ 2\sigma, \gamma\}$ and with $C$ depending only on $n$ and $\sigma$.
\end{teo}

\noindent
{\bf Proof.} For simplicity, we will see the case $2\sigma < 1$ and $2\sigma + \gamma < 1$. Defining $w = u_\epsilon - u$, for $x \in \Omega$ we have
\begin{eqnarray*}
- \I_\epsilon[w](x) &=& \I_\epsilon[u](x) + (-\Delta)^\sigma[u](x) \\
&=& -\epsilon^{n + 2\sigma} \int_{\Omega - x} \frac{u(x + z) - u(x)}{|z|^{n + 2\sigma} (\epsilon^{n+2\sigma} + |z|^{n + 2\sigma})}dz \\
&& -\epsilon^{n + 2\sigma} u(x) \int_{(\Omega - x)^c} \frac{dz}{|z|^{n + 2\sigma} (\epsilon^{n+2\sigma} + |z|^{n + 2\sigma})} \\
&=& I_1 + I_2.
\end{eqnarray*}
By the regularity of $u$ we have
\begin{eqnarray*}
|I_1| \leq C { \|u\|_{C^{2\sigma + \gamma}(\bar{\Omega})}} \epsilon^{\gamma},
\end{eqnarray*}
where $C$ does not depend on $\epsilon$. On the other hand,  for $I_2$ we split the analysis. First,  if $\epsilon \leq d(x)$, then
\begin{eqnarray*}
|I_2| \leq C{ \|u\|_{C^{2\sigma + \gamma}(\bar{\Omega})}} \epsilon^{n + 2\sigma}d(x)^{-(n + 2\sigma) + \gamma},
\end{eqnarray*}
where we have used that there is no loss of boundary condition for $u$, hence $u = 0$ on $\partial \Omega$ and then
$|u(x)| \leq [u]_{C^{2\sigma + \gamma}(\bar{\Omega})}d(x)^{2\sigma + \gamma}$. Hence, we conclude
$$
|I_2| \leq C \epsilon^{\gamma}.
$$
Second, when $d(x) < \epsilon$ we have
\begin{eqnarray*}
|I_2| \leq C d(x)^{2\sigma + \gamma} \epsilon^{-2\sigma} (d(x)/\epsilon)^{-2\sigma} \leq C\epsilon^\gamma.
\end{eqnarray*}
Since we know that $|w|\leq C\epsilon^{\beta_0}$ on $\partial \Omega$, by Proposition~\ref{propaprioriborde}, we can get the result proceeding 
exactly as in the proof of Proposition~\ref{cotauniformesolucion}.
\qed

\subsection{An example of a scheme without boundary equicontinuity.}

In this subsection we  consider the reverse scheme,
that is approximating zero order equations by fractional ones and we prove the absence of uniform modulus of continuity in $\bar{\Omega}$.  
For this, we recall some facts of Section \S \ref{ppiosdelmaximo}. Let $f \in C(\bar{\Omega})$ with $f \geq \varrho_0 > 0$, 
$J: \R^N \to \R_+$ integrable and $\I_J$ as in~\eqref{IJ}. Consider the associated problem~\eqref{eqbasica}-\eqref{exteriordata}, that is
\begin{eqnarray}\label{zeroorder}
\left \{ \begin{array}{rll} - \I_J[u] & = f \quad & \mbox{in} \ \Omega \\ u & = 0. \quad & \mbox{in} \ \Omega^c \end{array} \right .
\end{eqnarray}

As we saw in Remark~\ref{rmku>0borde}, the unique 
solution $u \in C(\bar{\Omega})$ for this problem is such that $u > 0$ in $\partial \Omega$.

Consider $J, f$ as above, with $J$ such that $J \geq m$ in $B_r$, for some $r, m > 0$. For $\epsilon \in (0,1)$ 
and $\alpha > 1$ consider the family of kernels
\begin{eqnarray*}
J_\epsilon(z) = \min \{1, |z/\epsilon|^\alpha\}^{-1}J(z),
\end{eqnarray*}
which are not integrable at the origin.  If we define
$$
\mathcal{J}_\epsilon[u](x) = \int_{\R^N} [u(x + z) - u(x)]J_\epsilon(z)dz 
$$
and consider the problems
\begin{eqnarray}\label{fractozero}
\left \{ \begin{array}{rll} -\mathcal{J}_\epsilon[u] & = f \quad & \mbox{in} \ \Omega \\ u & = 0, \quad & \mbox{in} \ \Omega^c \end{array} \right .
\end{eqnarray}
it is known that the
unique viscosity solution $u_\epsilon$ of~\eqref{fractozero} agrees the prescribed value of the equation on the boundary, and then 
$u_\epsilon = 0$ on $\partial \Omega$ for all $\epsilon \in (0,1)$, see for example \cite{Barles-Chasseigne-Imbert}. We have $\{ u_\epsilon \}$  
is uniformly bounded in $L^\infty(\bar{\Omega})$ and therefore, 
the application of half-relaxed limits together with
viscosity stability results in~\cite{Barles-Imbert}, imply $u_\epsilon \to u$ locally uniform in $\Omega$ as $\epsilon \to 0$, 
where $u$ is the unique solution to~\eqref{zeroorder}. Since $u$ is strictly positive in $\partial \Omega$, the convergence
of $u_\epsilon$ to $u$ cannot be uniform in $\bar{\Omega}$, and therefore the family $\{ u_\epsilon \}$ is not equicontinuous in this case.

This example resembles the behavior of the viscosity solutions $u_\epsilon$ of the equation 
$$
-\epsilon u'' + u' = 1 \quad \mbox{in} \ (0,1), \quad \mbox{with} \ u(0) = u(1) = 0, 
$$
which approximate the solution of the equation
$$
u' = 1 \quad \mbox{in} \ (0,1), \quad \mbox{with} \ u(0) = u(1) = 0.
$$
In this case, the family $(u_\epsilon)$ is not equicontinuous too, see~\cite{Barles-book}.

\bigskip

\noindent {\bf Acknowledgements:} The authors want to thank the referees for carefully reading the paper and making suggestions  
that resulted in a great improvement in the clarity and simplicity of the proofs. 
P.F.  was  partially supported by Fondecyt Grant \# 1110291,
  BASAL-CMM projects and CAPDE, Anillo ACT-125.
E.T. was partially supported by CONICYT, Grants Capital Humano Avanzado and Ayuda Realizaci\'on Tesis Doctoral.

\end{document}